\documentclass[12pt,a4paper]{article}
\usepackage{amssymb}
\usepackage{amsfonts}
\usepackage{amsmath}
\usepackage{hyperref}
\usepackage{geometry}

\input{tcilatex}
\begin{document}

\title{\textbf{Nystr\"{o}m Methods in the RKQ Algorithm for Initial-Value
Problems}}
\author{J.S.C. Prentice \\
%EndAName
Senior Research Officer\\
Mathsophical Ltd.\\
Johannesburg, South Africa}
\maketitle

\begin{abstract}
We incorporate explicit Nystr\"{o}m methods into the RKQ algorithm for
stepwise global error control in numerical solutions of initial-value
problems. The initial-value problem is transformed into an explicitly
second-order problem, so as to be suitable for Nystr\"{o}m integration. The
Nystr\"{o}m methods used are fourth-order, fifth-order and 10th-order. Two
examples demonstrate the effectiveness of the algorithm.
\end{abstract}

\section{Introduction}

In two previous papers we have considered the RK$rv$Q$z$ algorithm for
stepwise control of the global error in the numerical solution of an
initial-value problem (IVP), using Runge-Kutta methods \cite{prentice 2}\cite%
{prentice 3}. In the current paper, the third in the series, we focus our
attention on the use of Nystr\"{o}m methods in this error control algorithm
for $n$-dimensional problems of the form%
\begin{align}
\mathbf{y}^{\prime \prime }\left( x\right) & =\mathbf{f}\left( x,\mathbf{y}%
\right)   \label{ddy=f} \\
\mathbf{y}\left( x_{0}\right) & =\mathbf{y}_{0}  \notag \\
\mathbf{y}^{\prime }\left( x_{0}\right) & =\mathbf{y}_{0}^{\prime }.  \notag
\end{align}%
Note that $\mathbf{f}$ is not explicitly dependent on $\mathbf{y}^{\prime }$
(we note that Nystr\"{o}m methods can be used to solve the more general
problem $\mathbf{y}^{\prime \prime }\left( x\right) =\mathbf{f}\left( x,%
\mathbf{y,y}^{\prime }\right) ,$ but that is not our focus here). We
designate this Nystr\"{o}m-based algorithm RKN$rv$Q$z,$ and we will show in
a later section how any\ first-order IVP can be written in the form (\ref%
{ddy=f}), so that RKN$rv$Q$z$ is, in fact, generally applicable. The
motivation for considering this modification to RK$rv$Q$z$ is twofold: most
physical systems are described by second-order differential equations, and
Nystr\"{o}m methods applied to (\ref{ddy=f}) tend to be more efficient than
their Runge-Kutta counterparts.

\section{Relevant Concepts, Terminology and Notation}

Here we describe concepts, terminology and notation relevant to our work.
Note that boldface quantities are $n\times 1$ vectors, except for $\mathbf{%
\alpha }_{i}^{r},\mathbf{I}_{n},\mathbf{F}_{y}^{r},\mathbf{F}_{y^{\prime
}}^{r}$ and $\mathbf{g}_{y},$ which are $n\times n$ matrices.

\subsection{Nystr\"{o}m Methods \protect\cite{nystrom 1}\protect\cite{hairer}%
}

The most general definition of a Nystr\"{o}m method (sometimes known as
Runge-Kutta-Nystr\"{o}m (RKN)) for solving (\ref{ddy=f}) is%
\begin{align}
\mathbf{k}_{p}& =\mathbf{f}\left( x_{i}+c_{p}h_{i},\mathbf{w}_{i}+c_{p}h_{i}%
\mathbf{w}_{i}^{\prime }+h_{i}^{2}\sum\limits_{q=1}^{m}a_{pq}\mathbf{k}%
_{q}\right) \text{ \ \ \ \ \ \ \ }p=1,2,...,m  \notag \\
\mathbf{w}_{i+1}& =\mathbf{w}_{i}+h_{i}\mathbf{w}_{i}^{\prime
}+h_{i}^{2}\sum\limits_{p=1}^{m}b_{p}\mathbf{k}_{p}\equiv \mathbf{w}%
_{i}+h_{i}\mathbf{F}\left( x_{i},\mathbf{w}_{i}\right)  \label{RK defn 1} \\
\mathbf{w}_{i+1}^{\prime }& =\mathbf{w}_{i}^{\prime
}+h_{i}\sum\limits_{p=1}^{m}\widehat{b}_{p}\mathbf{k}_{p}.  \notag
\end{align}%
The coefficients $c_{p},a_{pq},b_{p}$ and $\widehat{b}_{p}$ are unique to
the given method. If $a_{pq}=0$ for all $p\leqslant q,$ then the method is
said to be \textit{explicit}; otherwise, it is known as an \textit{implicit}
RKN method. We will focus our attention on explicit methods. In the second
line of (\ref{RK defn 1}), we have implicitly defined the function $\mathbf{F%
}$. We treat $\mathbf{w}_{i}^{\prime }$ as an `internal parameter'; for our
purposes here, we do not identify $\mathbf{w}^{\prime }$ with $\mathbf{y}%
^{\prime },$ because $\mathbf{f}$ is not dependent on $\mathbf{y}^{\prime }$%
. The symbol $\mathbf{w}$ is used here and throughout to indicate the
approximate numerical solution, whereas the symbol $\mathbf{y}$ will be used
to denote the exact solution. We will denote an RKN method of order $r$ as
RKN$r$ and, for such a method, we write%
\begin{equation}
\mathbf{w}_{i+1}^{r}=\mathbf{w}_{i}^{r}+h_{i}\mathbf{F}^{r}\left( x_{i},%
\mathbf{w}_{i}^{r},\mathbf{w}_{i}^{r\prime }\right) .  \label{w^r = ...}
\end{equation}%
The stepsize $h_{i}$ is given by%
\begin{equation*}
h_{i}\equiv x_{i+1}-x_{i}
\end{equation*}%
and carries the subscript because it may vary from step to step. It is known
that RKN$r$ has a local error of order $r+1$ and a global error of order $r,$
just like its Runge-Kutta counterpart RK$r$.

\subsection{IVPs in the form $y^{\prime \prime }=f\left( x,\mathbf{y}\right) 
$}

Consider the $n$-dimensional IVP%
\begin{align}
\mathbf{y}^{\prime }\left( x\right) & =\mathbf{g}\left( x,\mathbf{y}\right) 
\label{dy=g} \\
\mathbf{y}\left( x_{0}\right) & =\mathbf{y}_{0}.  \notag
\end{align}%
This gives%
\begin{equation*}
y_{j}^{\prime \prime }=\sum\limits_{i=1}^{n}\frac{\partial g_{j}\left( x,%
\mathbf{y}\right) }{\partial x}+\frac{\partial g_{j}\left( x,\mathbf{y}%
\right) }{\partial y_{i}}\frac{dy_{i}}{dx}=\sum\limits_{i=1}^{n}\frac{%
\partial g_{j}\left( x,\mathbf{y}\right) }{\partial x}+\frac{\partial
g_{j}\left( x,\mathbf{y}\right) }{\partial y_{i}}g_{i}\left( x,\mathbf{y}%
\right) 
\end{equation*}%
where $y_{i}$ is the $i$th component of $\mathbf{y,}$ and $g_{i}$ is the $i$%
th component of $\mathbf{g.}$ Clearly, we have%
\begin{equation*}
y_{j}^{\prime \prime }=\sum\limits_{i=1}^{n}\frac{\partial g_{j}\left( x,%
\mathbf{y}\right) }{\partial x}+\frac{\partial g_{j}\left( x,\mathbf{y}%
\right) }{\partial y_{i}}g_{i}\left( x,\mathbf{y}\right) \equiv f_{j}\left(
x,\mathbf{y}\right) 
\end{equation*}%
for all $j=1,2,\ldots ,n,$ and so we can write%
\begin{equation*}
\mathbf{y}^{\prime \prime }\left( x\right) =\mathbf{f}\left( x,\mathbf{y}%
\right) .
\end{equation*}%
The initial values for this second-order problem are then given by%
\begin{align*}
\mathbf{y}\left( x_{0}\right) & =\mathbf{y}_{0} \\
\mathbf{y}^{\prime }\left( x_{0}\right) & =\mathbf{g}\left( x_{0},\mathbf{y}%
_{0}\right) \equiv \mathbf{y}_{0}^{\prime }.
\end{align*}%
Hence, any first-order IVP can be transformed into an IVP of the form (\ref%
{ddy=f}). This is ideally suited to the Nystr\"{o}m methods, which are
specifically designed for this type of IVP. They are also more efficient
than their Runge-Kutta counterparts; for example, the methods to be used
later, RKN4 and RKN5, require three and four stage evaluations,
respectively, as opposed to RK4 and RK5, which require at least four and six
stage evaluations, respectively.

\subsection{Error Propagation in RKN}

It can be shown \cite{prentice 1} that, for RK$r$,%
\begin{align}
\mathbf{\Delta }_{i+1}^{r}& \equiv \mathbf{w}_{i+1}^{r}-\mathbf{y}_{i+1}=%
\mathbf{\varepsilon }_{i+1}^{r}+\mathbf{\alpha }_{i}^{r}\mathbf{\Delta }%
_{i}^{r}  \label{master error equn} \\
\mathbf{\alpha }_{i}^{r}& \equiv \mathbf{I}_{n}+h_{i}\mathbf{F}%
_{y}^{r}\left( x_{i},\mathbf{\xi }_{i}\right) ,  \label{alpha r =}
\end{align}%
where $\mathbf{\varepsilon }_{i+1}^{r}=O\left( h_{i}^{r+1}\right) $ is the
local error, $\mathbf{\Delta }_{i+1}^{r}$ is the global error and $\mathbf{F}%
_{y}^{r}$ is the Jacobian (with respect to $\mathbf{y})$ of the function $%
\mathbf{F}^{r}\left( x_{i},\mathbf{w}_{i}^{r}\right) $ associated with RK$r$%
. The term $h_{i}\mathbf{F}_{y}^{r}\left( x_{i},\mathbf{\xi }_{i}\right) $
in the matrix $\mathbf{\alpha }_{i}^{r}$ arises from a first-order Taylor
expansion of $\mathbf{F}^{r}\left( x_{i},\mathbf{w}_{i}\right) =$ $\mathbf{F}%
^{r}\left( x_{i},\mathbf{y}_{i}+\mathbf{\Delta }_{i}^{r}\right) $ with
respect to $\mathbf{y}_{i}$.

For a Nystr\"{o}m method RKN$r,$ we have $\mathbf{F}^{r}=\mathbf{F}%
^{r}\left( x_{i},\mathbf{w}_{i}^{r}\right) $ and so, as above,%
\begin{equation*}
\mathbf{\alpha }_{i}^{r}\equiv \mathbf{I}_{n}+h_{i}\mathbf{F}_{y}^{r}\left(
x_{i},\mathbf{\zeta }_{i}\right) ,
\end{equation*}%
where $\mathbf{\zeta }_{i}$\textbf{\ }is an appropriate constant. Hence, the
global error in RKN$r$ is also given by (\ref{master error equn}).

\subsection{RK$rv$Q$z$}

We will not discuss RK$rv$Q$z$ in detail here; the reader is referred to our
previous work where the algorithm has been discussed extensively. It
suffices to say that RK$rv$Q$z$ uses RK$r$ and RK$v$ to control local error
via local extrapolation, while simultaneously using RK$z$ to keep track of
the global error in the RK$r$ solution. Such global error arises due to the
propagation of the RK$v$ global error. RK$rv$Q$z$ is designed to estimate
the various components of the global error in RK$r$ and RK$v$ at each node
and, when the global error is deemed too large, a quenching procedure is
carried out. This simply involves replacing the RK$r$ and RK$v$ solutions
with the much more accurate RK$z$ solution, whenever necessary, so that the
RK$r$ and RK$v$ global errors do not accumulate beyond a desired tolerance.

\subsection{RKN$rv$Q$z$}

The algorithm RKN$rv$Q$z$ is nothing more than RK$rv$Q$z$ with RK$r$, RK$v$
and RK$z$ replaced with RKN$r$, RKN$v$ and RKN$z$. Of course, RKN$rv$Q$z$ is
applied to problems of the form (\ref{ddy=f}), whereas RK$rv$Q$z$ is applied
to problems of the form (\ref{dy=g}).

We also report on a refinement to the algorithm: in RK$rv$Q$z,$ if the
global error at $x_{i}$ is too large, we replace $\mathbf{w}_{i}^{r}$ with $%
\mathbf{w}_{i}^{z}$ and then recompute $\mathbf{w}_{i+1}^{r}$ and $\mathbf{w}%
_{i+1}^{v},$ using $\mathbf{w}_{i}^{z}$ as input for both RK$r$ and RK$v.$
This is the essence of the quenching procedure. However, in retrospect it
seems quite acceptable to simply replace $\mathbf{w}_{i+1}^{r}$ and $\mathbf{%
w}_{i+1}^{v}$ with $\mathbf{w}_{i+1}^{z};$ this avoids the need for
recomputing $\mathbf{w}_{i+1}^{r}$ and $\mathbf{w}_{i+1}^{v},$ which
improves efficiency and, after all, it is the global error in $\mathbf{w}%
_{i+1}^{r}$ and $\mathbf{w}_{i+1}^{v},$ not $\mathbf{w}_{i}^{r}$ and $%
\mathbf{w}_{i}^{v},$ that is too large. Both approaches are effective,
although one is more efficient than the other. It is the more efficient
approach that we have employed in RKN$rv$Q$z.$

\section{Numerical Examples}

It is not our intention to compare methods or algorithms but, for the sake
of consistency, we will apply RKN$rv$Q$z$ to the same examples that we
considered in our previous work on RK$rv$Q$z.$ In our calculations, we use
RKN4, RKN5 and RKN10 which gives the algorithm RKN45Q10. RKN4 and RKN5 are
taken from Hairer et al\textit{\ }\cite{hairer}, and RKN10 is from Dormand
et al \cite{dormand}.

The first of these is the scalar problem%
\begin{align*}
y^{\prime }& =\left( \frac{\ln 1000}{100}\right) y \\
y\left( 0\right) & =1
\end{align*}%
which transforms to%
\begin{align*}
y^{\prime \prime }& =\left( \frac{\ln 1000}{100}\right) ^{2}y \\
y\left( 0\right) & =1 \\
y^{\prime }\left( 0\right) & =\frac{\ln 1000}{100}.
\end{align*}%
Solving this problem with RKN45 and RKN45Q10 with a tolerance of $10^{-10}$
on the absolute local and global errors gives the error curves shown in
Figure 1. The global error obtained with RKN45 is clearly larger than the
desired tolerance on most of the interval, despite local error control via
local extrapolation. However, RKN45Q10 yields a solution with a global error
always less than the tolerance - the maximum global error in this case is $%
9.1\times 10^{-11}$. The points on the $x$-axis where this global error
decreases sharply correspond to the quenches carried out using RKN10.

The second example is the simple harmonic oscillator%
\begin{align*}
y_{1}^{\prime }& =y_{2} \\
y_{2}^{\prime }& =-y_{1} \\
\mathbf{y}\left( 0\right) & =\left[ 
\begin{array}{c}
0 \\ 
1000%
\end{array}%
\right]
\end{align*}%
which has solution%
\begin{eqnarray*}
y_{1}\left( x\right) &=&1000\sin x \\
y_{2}\left( x\right) &=&1000\cos x
\end{eqnarray*}%
and becomes, in explicit second-order form,%
\begin{align*}
\mathbf{y}^{\prime \prime }& =\left[ 
\begin{array}{c}
y_{1}^{\prime \prime } \\ 
y_{2}^{\prime \prime }%
\end{array}%
\right] =\left[ 
\begin{array}{c}
-y_{1} \\ 
-y_{2}%
\end{array}%
\right] \equiv \mathbf{f}\left( x,\mathbf{y}\right) \\
\mathbf{y}\left( 0\right) & =\left[ 
\begin{array}{c}
0 \\ 
1000%
\end{array}%
\right] ,\text{ \ }\mathbf{y}^{\prime }\left( 0\right) =\left[ 
\begin{array}{c}
1000 \\ 
0%
\end{array}%
\right] .
\end{align*}%
Since the solution oscillates between $-1000$ and $1000$, there are regions
where the solution has magnitude less than unity - here, we implement
absolute error control - and regions where the solution has magnitude
greater than unity, where we implement relative error control. With an
imposed tolerance of $10^{-8}$ on the local and global errors (relative and
absolute) we found a maximum global error of $\sim 4\times 10^{-8}$ in each
component when using RKN45, and a global error no greater than $0.99\times
10^{-8}$ with RKN45Q10, on $x\in \left[ 0,200\right] .$ A total of 20
quenches were needed.

\section{Conclusion}

We have considered the use of Nystr\"{o}m methods in RK$rv$Q$z,$ wherein a
combination of local extrapolation and quenching result in stepwise global
error control in numerical solutions of IVPs. Two examples have demonstrated
the success of RKN45Q10.

\end{document}